%% file: super_1_latest.tex
\documentclass{ajmart}

\include{september}  \pagespan{115}{137}

\usepackage{amsmath,amssymb,amsbsy,amsfonts,amsthm,latexsym,amsopn,amstext, amsxtra,euscript,amscd} %, hyperref}

\usepackage{cite}

\usepackage{tikz}
\usetikzlibrary{positioning,shadows,arrows}

\usepackage{amsrefs}

\usepackage{xy} 
\input xy \xyoption{all}

\newtheorem{theorem}{Theorem}

\newtheorem{lemma}{Lemma}
\newtheorem{remark}{Remark}

\newtheorem{exa}{Example}
\theoremstyle{definition}

\def\diag{\mbox{diag }}

\def\sem{{\rtimes}}

\begin{document}

\title{On superelliptic curves of level $n$ and their quotients, I.}

\subjclass[2000]{14H32, 14H37, 14K25}

\keywords{cyclic quotients, algebraic curves, automorphism groups.}

%\hline 

\noindent \rule[-0.1cm]{13cm}{0.02cm} %Text\rule[-0.1cm]{3cm}{0.01cm}   

\vspace{1.6cm}

\maketitle

%, V. Hoxha,  T. Shaska
%
\begin{small}
\begin{center} 
\textsc{Lubjana Beshaj}\\
Department of Mathematics and Statistics,\\
Oakland University, Rochester, MI, 48309.\\
Email:  beshaj@oakland.edu \\ 
\end{center}

\medskip

\begin{center} 
\textsc{Valmira Hoxha}\\
Department of Mathematics and Statistics,\\
Oakland University, Rochester, MI, 48309.\\
Email:  vhoxhaj@oakland.edu \\
\end{center}

\medskip

\begin{center} 
\textsc{Tony Shaska}\\
Department of Mathematics and Statistics,\\
Oakland University, Rochester, MI, 48309.\\
Email:  shaska@oakland.edu \\
\end{center}

\end{small}

\vspace{.5cm}

%\hline

\rule[-0.1cm]{11.5cm}{0.01cm}

\vspace{.2cm}

%****************
\begin{abstract} 
We study families of superelliptic curves with fixed automorphism groups.  Such families are parametrized with  invariants expressed in terms of the coefficients of the curves.  Algebraic relations among such invariants determine the lattice of inclusions among the loci of superelliptic curves and their field of moduli.  We give a Maple package of how to compute the normal form of an superelliptic curve and its  invariants.  A complete list of all superelliptic curves of genus $g \leq 10$ defined over any field of characteristic $\neq 2$  is given in a subsequent paper \cite{super2}. 
\end{abstract}

\vspace{0.2cm}  

%\hline

\rule[-0.1cm]{11.5cm}{0.01cm}

\vspace{1cm}
%\tableofcontents

%*******************  our definitions 

\def\P{\mathbb P^1 (k)}
\def\Q{{\mathbb Q}}
\def\Z{{\mathbb Z}}
\def\C{{\mathbb C}}
\def\M{{\mathcal M}}
\def\H{\mathcal H}
\def\L{\mathcal L}

\def\J{\mbox{Jac }}
\def\O{\Omega}
\def\T{\theta}

\def\X{\mathcal X}
\def\Y{\mathcal Y}
\def\A{\mathcal A}
\def\B{\mathcal B}

\def\Aut{\mbox{Aut }}
\def\bAut{\overline {\mathrm{Aut}}}

\def\Pic{\mbox{Pic }}
\def\Jac{\mbox{Jac }}

\def\iso{\equiv}

\def\embd{\hookrightarrow}

\newcommand\G{\bar G}
\def\D{\Delta}
\def\e{\varepsilon}
\def\u{\mathfrak s}

\def\g{\gamma}
\def\bG{\bar G}
\def\t{\tau}
\def\d{{\delta }}
\def\a{{\alpha }}
\def\b{{\beta }}
\def\c{\textbf{c}}

\def\bC{\mathfrak \si}

\def\th{\theta}
\def\<{\langle}
\def\>{\rangle}

\def\diag{\mbox{diag }}
\def\sem{{\rtimes}}

\def\s{\mathfrak s}

\def\v{\mathfrak v}

\def\r{\gamma}

\def\ss{$\mathfrak s$}

\def\p{\mathfrak p}

\newcommand{\cH}{\mathcal H}

\newcommand\F{\mathbb F}

\def\emb{\hookrightarrow }
\def\tr{{\mbox tr }}

\def\U{\mathcal U}

\def\normal{\triangleleft}

\def\xs{{\mathbb o}}
\def\xs{{\Bbb o}}

%\setpage{115}

%\medskip

%\hline

%\bigskip
%\tableofcontents

%&&&&&&&&&&&&&&&&&&&&&&&&&&&&&&&&&&&&&&&&&&&&&&&&&&&&&&&&&&&&&&&&

%antoniadis, kontogiorgis, dihedral-invariants, sanjeewa

\section{Introduction}

Let $\X_g$ be an  algebraic curve of genus $g \geq 2$ defined over an algebraically closed field $k$ of characteristic  $p \neq 2$. 
What is the group of automorphisms of $\X_g$ over $k$?  Given the group of automorphisms $G$ of a genus $g$ curve, can we determine the equation of the curve?  
These two questions have been studied for a long time and a complete answer is not known for either one.  There are some families of curves where we can answer completely the above questions, such as the hyperelliptic curves.   The Klein's curve was the first celebrated example of a non-hyperelliptic curve where the the automorphisms of the curve and its equation are shown explicitly.  The main purpose of this paper is to show that we can do this for a larger family of curves. 

In characteristic zero, the first question is answered by work of Magaard, Shaska, Shpectorov, V\"olklein, et al.  Based on previous work of Breuer and using computer algebra systems as GAP,  they show how one can compute the list of full automorphism groups for any fixed genus $g \geq 2$.  It is still an unsettled question the case of positive characteristic, where many  tedious cases of wild ramifications need to be considered.  The second question is unsettled even in characteristic zero.  It is much harder to determine a parametric equation for the curve, given its group of automorphisms $G$.  

However, if we go through the lists of groups $G$ which occur as automorphism groups of  genus $g$ curves we notice, as to be expected, that the majority of them have the following property;  there is a central element $\tau \in G$ such that    the quotient space $\X_g /\< \tau \>$ has genus zero. Such curves in the literature are called \textit{superelliptic curves} or \textit{cyclic curves}.  For the purposes of this paper we will use the term \textbf{superelliptic curves of level $n$}.

Hence, for a fixed genus $g$, certain families of curves have equation $y^n = f(x)$, for some $n\in \Z$ and a generic polynomial $f(x) \in k[x]$.  The values of $n$ depend solely on the genus $g$ and the field $k$.  Such cases we call them \textit{root cases} or \textit{fundamental cases}.  For a given $n$ let $\H_n$ denote a connected component  of the space of genus $g$ curves with equation as above.  Isomorphism classes of curves in $\H_n$ are determined by the invariants of degree $n$ binary forms.  Such invariants were the main focus of classical invariant theory in the 19-th century and they are only known for $n \leq 8$. Even for $n \leq 8$ the expressions of such invariants in terms of the coefficients of $f(x)$ are quite long and not so convenient for computations. 

If the curve has an additional automorphism then this automorphism has to permute the roots of $f(x)$.  In this case, additional invariants can be defined in terms of the coefficients of $f(x)$.  These invariants were first discovered by Shaska for genus two curves in \cite{sh_2000} and then generalized by Shaska/Gutierrez for all hyperelliptic curves in \cite{g_sh},  where they were called \textbf{dihedral invariants}.  Moreover, in \cite{g_sh} was determined a relation among such invariants, for any genus $g$, in the case of hyperelliptic curves with an extra involution.  
In \cite{GSS} algebraic relations among such invariants were computed for the case of  genus three hyperelliptic curves and a method was described how to compute such relations in general. Extending work done by Gutierrez/Shaska in   \cite{g_sh},  Antoniadis and Kontogeorgis defined these invariants   \cite{AK} for cyclic covers of $\mathbb P^1 (k)$ for  positive characteristic. In these paper we will call them \textbf{\ss-invariants} and will describe how to compute them for any genus $g$ superelliptic curve of level $n$. 

Superelliptic curves are quite important in many applications. They are the only curves where we fully understand the automorphism groups for every characteristic and can associate an equation of the curve in each case of the group.   The full groups of automorphisms of superelliptic curves defined over a field of characteristic zero followed from previous work of Magaard et al, \cite{kyoto}.  However, for the first time a complete list of full automorphism groups of superelliptic curves for odd characteristic was determined  by Sanjeewa in \cite{Sa1}.  The equations for each family, when the full automorphism group was fixed, were determined by Sanjeewa and Shaska in \cite{Sa2}.  Such curves were further studied in  \cite{nato_beshaj, beshaj}, where singular subloci of $\M_2$ were studied,  in applications in coding in \cite{elezi_sh}.  The \ss-invariants which we study in section 4 were discovered in \cite{sh_2000} and used by several authors in many applications since them.  For further applications of such invariants one can check   \cite{sevilla,   sh_2000, nato_wijesiri, serdica, sh_04, sh_03, issac, ajm_sh1, sh_05, deg3, serdica, open_problems,  sh_02, sh_01, super2, super3, super4}.

In this paper we give  a list of automorphism groups of superelliptic curves of genus $g$ and the corresponding equation for each group. We define invariants for such curves and give algorithms how to compute such invariants and how to determine algebraic relations among them.  Such computations are completely done in the case of genus 3, in order to provide some general idea of the genus $g > 3$ case.  

\bigskip

\noindent \textbf{Notation:} Throughout this paper by $g$ we denote an integer $\geq 2$ and  $k$ denotes an algebraically closed field of characteristic $\neq 2$.  Unless otherwise noted,  by a "curve" we always mean the isomorphism class of an algebraic curve defined over $k$.   The automorphism group of a curve always means the full automorphism group of the curve. 

%********************************************************************
\section{Preliminaries on automorphisms of the projective line.}

In this section we set the notation and describe briefly some general facts.  Fix an integer $g\ge2$.  Let $\X_g$ denote a genus $g$ generic curve defined over an algebraically closed field $k$ of characteristic $p \geq 0$. We denote by $G$ the full automorphism group of $\X_g$.  Hence,  $G$ is a finite group. Denote by $K$ the function field of $\X_g$ and assume that the affine equation of $\X_g$ is given some polynomial in terms of $x$ and $y$. 

Let $H=\< \tau \>$ be a cyclic subgroup of $G$ such that $| H | = n$ and $H \normal G$, where $n \geq 2$. Moreover, we assume that the quotient curve $\X_g / H$ has genus zero. 
The \textbf{reduced automorphism group of $\X_g$ with respect to $H$} is called the group  $\G \, := \, G/H$, see \cite{K, Sa1}.  

Assume  $k(x)$ is the  genus zero subfield of $K$ fixed by $H$.   Hence, $[ K : k(x)]=n$. Then, the group  $\G$ is a subgroup of the group of automorphisms of a genus zero field.   Hence, $\G <  PGL_2(k)$ and $\G$ is finite. It is a classical result that every finite subgroup of $PGL_2 (k)$  is  isomorphic to one of the following: $C_m $, $ D_m $, $A_4$, $S_4$, $A_5$, \emph{semidirect product  of elementary  
Abelian  group  with  cyclic  group}, $ PSL(2,q)$ and $PGL(2,q)$,  see \cite{VM}.

The group $\G$ acts on $k(x)$ via the natural way. The fixed field of this action is a genus 0 field, say $k(z)$. Thus, $z$
is a degree $|\G| := m$ rational function in $x$, say $z=\phi(x)$. We illustrate with the following diagram:

\begin{figure}[ht]
\[
\xymatrix{
K =k(x, y)\ar@{-}[d]_{\, \, n }^{\, \, \, H} \\
k(x)=k(x, y^n)\ar@{-}[d]_{\, \, m }^{\, \, \, \G}  \\ 
E=k (z) \\ 
}
\qquad \qquad 
\xymatrix{  {\X_g}\ar[d]_{\, \, \phi_0 }^{\, \, \,  H} \\ 
\P\ar[d]_{\, \, \, \phi}^{\, \, \, \G}  \\ 
\P \\
} 
\]
\caption{The automorphism groups and the corresponding covers}
\label{fig1}
\end{figure}

It obvious that $G$ is a degree $n$ extension of $\G$ and $\G$ is a finite subgroup of $PGL_2(k)$.  Hence, if we know all the possible groups that occur as $\G$ then hopefully we can figure out $G$ and the equation for $K$. 

To do this we have to recall some classical result on finite subgroups of the projective linear group $PGK_2(k)$ and their fixed fields.   
First we define a semidirect product of an elementary Abelian group  with  a cyclic group as follows, see \cite{VM} for details. 

Let $\mbox{char } k=p$ and $k = \F_q$ for $q= p^r$.       For each $m \, | \left( p^t - 1 \right)$, $t = 1, \dots , r$,   we define $\U_m$ as follows 
\[\U_m :=\{a \in k \, | \, (a \prod_{j=0}^{\frac{p^t-1}{m}-1}(a^m-b_j))=0, \, b_j \in  k^* \}. \]
Obviously $\U_m$ is a subgroup of the additive group of $k$.    Let
\[K_m :=\<  \{ \sigma_a(x)=x+a, \, \tau(x)=\xi^2x    \,     |  \,  \forall    a \in \U_m  \}  \> , \]
where   $\xi$ is a primitive $2m$-th root of unity. 
Now we are ready to state the following classical result. 

\begin{theorem}\label{thm_1}
i) Let $k$ be an algebraically closed field of characteristic $p\neq 2$ of size $q$ when $k$ is finite and $G$ be a finite subgroup of $PGL_2(k)$. Then, $G$ is isomorphic to one of the following groups 
\[C_m, D_m, A_4, S_4, A_5, U = C_p^t, K_m, PSL_2(q), \textit{ and } PGL_2(q),\]
where $(m, p)=1$ and $C_m$ (resp. $D_m$) denotes the cyclic (resp. dihedral) group of size $m$. 

ii) Let $G$ act on $k(x)$ in the natural way.  The fixed field of $G$ is a genus zero subfield $k(z)$, where $z$  is given as in  Table ~\ref{t1}, with $\alpha =\frac{q(q-1)}{2}, \quad \beta= \frac{q+1}{2}$ and 
$H_t$ is a subgroup of the additive group of $k$ with order $| H_t | = p^t$ and $b_j \in k^*$.
\end{theorem}
\begin{table}[hb!]  
\begin{center}
\begin{tabular}{c|ccc}
$Case$ & $\G$ & $z=\phi (x)$ & $Ramification$  \\
\hline \\
1 & $C_m$, $(m,p)=1$& $x^m$ & $(m,m)$\\ \\
2 & $D_{2m}$, $(m,p)=1$& $x^m+\frac{1}{x^m}$ & $(2,2,m)$\\ \\
3 & $A_4, \, p\neq 2, 3$ & $\frac{x^{12}-33x^8-33x^4+1}{x^2(x^4-1)^2}$ & $(2,3,3)$\\ \\
4 & $S_4, \, p\neq 2, 3$ & $\frac{(x^8+14x^4+1)^3}{108(x(x^4-1))^4}$ & $(2,3,4)$\\ \\
5 & $A_5, \,p\neq 2, 3, 5$ & $\frac{(-x^{20}+228x^{15}-494x^{10}-228x^5-1)^3}{(x(x^{10}+11x^5-1))^5}$ & $(2,3,5)$\\ \\
  & $A_5, \,p=3$ & $\frac{(x^{10}-1)^6}{(x(x^{10}+2ix^5+1))^5}$ & $(6,5)$\\ \\
6 & $U$ & $  \displaystyle{\prod_{a \in H_t}} (x+a)$ & $(p^t)$\\ \\
7 & $K_m$ & $(x   \displaystyle{\prod_{j=0}^{\frac{p^t-1}{m}-1} } (x^m-b_j))^m$ & $(mp^t,m)$ \\ \\
8 & $PSL(2,q), \,p\neq 2$ & $\frac{((x^q-x)^{q-1}+1)^{\frac{q+1}{2}}}{(x^q-x)^{\frac{q(q-1)}{2}}}$ & $(\alpha,\beta)$ \\ \\
9 & $PGL(2,q)$ & $\frac{((x^q-x)^{q-1}+1)^{q+1}}{(x^q-x)^{q(q-1)}}$ & $(2\alpha,2\beta)$ \\ \\
\end{tabular}
\vspace{.5cm}
\caption{Rational functions for each finite  $G < PGL_2 (k)$} \label{t1}
\end{center}
\end{table}

Proof of the first part can be found in \cite{VM} and verifying the second part is an easy computational exercise.  Next, we continue with our tasks of determining $G$ and an equation for $K$. 

Let $\phi_0: \X_g \to \P$ be the cover which corresponds to the degree n extension $K/k(x)$. Then $\Phi :=\phi \circ \phi_0$ has monodromy group $G:=\Aut(\X_g)$. From the basic covering theory, the group $G$ is embedded in the group $S_l$ where $l=\deg \Phi = n m$. There is an $r$-tuple $\overline{\sigma}:=(\sigma_1, \dots ,\sigma_r)$, where $\sigma_i \in S_l$ such that $\sigma_1,...,\sigma_r$ generate $G$ and $\sigma_1, \dots ,\sigma_r=1$. The signature of $\Phi$ is an $r$-tuple of conjugacy classes 
\[\sigma:=(C_1, \dots ,C_r)\]
in $S_l$ such that $C_i$ is the conjugacy class of $\sigma_i$. 
We use the notation $i$ to denote the conjugacy class of permutations which is cycle of length $i$. 
Using the signature of $\phi: \P \to \P$ one finds out the signature of $\Phi :\X_g \to \P$ for any given $g$ and $G$.

For the extension $K / E$, from  the Hurwitz genus formula we have that
\begin{equation}\label{e1}
2(g_K-1)=2(g_E-1)|G|+deg(\mathfrak{D}_{K/E})
\end{equation}
with $g_K$ and $g_E$ the genera of $K$ and $E$ respectively and $\mathfrak{D}_{K/E}$ the different of $K/E$. Let $\overline{P}_1, \overline{P}_2, ... , \overline{P}_r$ be ramified primes of $E$. If we set $d_i=deg(\overline{P}_i)$ and let $e_i$ be the ramification index of the $\overline{P}_i$ and let $\beta_i$ be the exponent of $\overline{P}_i$ in $\mathfrak{D}_{K/E}$. Hence, (1) may be written as
\begin{equation}\label{e2}
2(g_K-1)=2(g_E-1)|G|+|G|\sum_{i=1}^{r}\frac{\beta_i}{e_i}d_i
\end{equation}

If $\overline{P}_i$ is tamely ramified then $\beta_i=e_i-1$ or if $\overline{P}_i$ is wildly ramified then $\beta_i=e_i^*q_i+q_i-2$ with $e_i=e_i^*q_i$, $e_i^*$ relatively prime to $p$, $q_i$ a power of $p$ and $e_i^*|q_i-1$. 

For a fixed $G$ and  $\sigma$ the family of covers $\Phi : \X_g\to \P$ is a Hurwitz space $\cH(G,\sigma)$. $\cH(G,\sigma)$ is an irreducible algebraic variety of dimension $\d(G,\sigma)$. Using Eq.~\eqref{e2} and signature $\sigma$ one can find out the dimension for each $G$.

%******************************************************************

\section{Superelliptic curves}

The superelliptic curves by definition have the group $H$ as a subgroup of their automorphism group. However, the curve might have more automorphisms.  Determining the full automorphism group is equivalent to determine degree $n$ extensions of $\G$, where $G$ is as above.  

The following  theorems give us all possible automorphism groups of genus $g\geq 2$ superelliptic curves defined over any $k$ such that $char k \neq 2$, see \cite{K, Sa1, Sa2} for details. 

\begin{theorem}[Sanjeewa, 2010] \label{thm_2}
Let $\X_g$ be a genus $g\geq2$ irreducible superelliptic curve defined over an algebraically closed field $k$, $\mbox{char } (k)=p\neq2$. Let  $G=Aut(\X_g)$, $\G$ its reduced automorphism group with respect to $H$, where $|H|=n$. Then, $G$ is isomorphic to one of the following:\\
\begin{enumerate}
\item If $\G \cong C_m$ then $G \cong C_{mn}$ or $G$ is isomorphic to 
\[\left\langle \r, \sigma \right|\r^n=1,\sigma^m=1,\sigma\r\sigma^{-1}=\r^l \rangle, \]
where $(l,n)=1 \textit{ and }   l^m\equiv 1 \mod n$. \\

\item If $\G \cong D_{2m}$ for some $m \in \Z$, then $G \cong D_{2m} \times C_n$, $G \cong D_{2mn}$, or   $G$ is isomorphic to 
\begin{align*}
\begin{split}
i) & \left\langle \r, \sigma, \tau \right|\r^n=1,\sigma^2=\r,\t^2=1,(\sigma\t)^m=1,\sigma\r\sigma^{-1}=\r,\t\r\t^{-1}=\r^{n-1} \rangle\\
ii) & \left\langle \r, \sigma, \tau \right|\r^n=1,\sigma^2=\r,\t^2=\r^{n-1},(\sigma\t)^m=1,\sigma\r\sigma^{-1}=\r,\t\r\t^{-1}=\r \rangle\\
iii) & \left\langle \r, \sigma, \tau \right|\r^n=1,\sigma^2=\r,\t^2=1,(\sigma\t)^m=\r^{\frac{n}{2}},\sigma\r\sigma^{-1}=\r,\t\r\t^{-1}=\r^{n-1} \rangle \\
iv) & \left\langle \r, \sigma, \tau \right|\r^n=1,\sigma^2=\r,\t^2=\r^{n-1},(\sigma \t)^m=
\r^{\frac{n}{2}},\sigma \r\sigma^{-1}=\r,\t\r\t^{-1}=\r \rangle \\
\end{split}
\end{align*} \\

\item If $\G \cong A_4$ and $p\neq3$ then $G \cong A_4 \times C_n$ or   $G$ is isomorphic to 
\begin{align*}
\begin{split}
i) & \left\langle \r, \sigma, \tau \right|\r^n=1,\sigma^2=1,\t^3=1,(\sigma\t)^3=1,\sigma\r\sigma^{-1}=\r,\t\r\t^{-1}=\r^l \rangle\\
ii) & \left\langle \r, \sigma, \tau \right|\r^n=1,\sigma^2=1,\t^3=\r^{\frac{n}{3}},(\sigma\t)^3=\r^{\frac{n}{3}},\sigma\r\sigma^{-1}=\r,\t\r\t^{-1}=\r^l \rangle\\
\end{split}
\end{align*}
where $(l,n)=1$ and $l^3\equiv 1 \mod n$ or
\begin{center}
$\left\langle \r, \sigma, \tau \right|\r^n=1,\sigma^2=\r^{\frac{n}{2}},\t^3=\r^{\frac{n}{2}},(\sigma\t)^5=\r^{\frac{n}{2}},\sigma\r\sigma^{-1}=\r,\t\r\t^{-1}=\r \rangle $
\end{center}
or
\begin{align*}
\begin{split}
iii) & \left\langle \r, \sigma, \tau \right|\r^n=1,\sigma^2=1,\t^3=1,(\sigma\t)^3=1,\sigma\r\sigma^{-1}=\r,\t\r\t^{-1}=\r^k \rangle\\
iv) & \left\langle \r, \sigma, \tau \right|\r^n=1,\sigma^2=\r^{\frac{n}{2}},\t^3=1,(\sigma\t)^3=1,\sigma\r\sigma^{-1}=\r,\t\r\t^{-1}=\r^k \rangle\\
\end{split}
\end{align*}
where $(k,n)=1$ and $k^3\equiv 1 \mod n$. \\

\item If $\G \cong S_4$ and $p\neq3$ then $G \cong S_4 \times C_n$ or  $G$ is isomorphic to 
\begin{align*}
\begin{split}
i) & \left\langle \r, \sigma, \tau \right|\r^n=1,\sigma^2=1,\t^3=1,(\sigma\t)^4=1,\sigma\r\sigma^{-1}=\r^l,\t\r\t^{-1}=\r \rangle\\
ii) & \left\langle \r, \sigma, \tau \right|\r^n=1,\sigma^2=1,\t^3=1,(\sigma\t)^4=\r^{\frac{n}{2}},\sigma\r\sigma^{-1}=\r^l,\t\r\t^{-1}=\r \rangle\\
iii) & \left\langle \r, \sigma, \tau \right|\r^n=1,\sigma^2=\r^{\frac{n}{2}},\t^3=1,(\sigma\t)^4=1,\sigma\r\sigma^{-1}=\r^l,\t\r\t^{-1}=\r \rangle \\
iv) & \left\langle \r, \sigma, \tau \right|\r^n=1,\sigma^2=\r^{\frac{n}{2}},\t^3=1,(\sigma\t)^4=\r^{\frac{n}{2}},\sigma\r\sigma^{-1}=\r^l,\t\r\t^{-1}=\r \rangle\\
\end{split}
\end{align*}
where $(l,n)=1$ and $l^2\equiv 1 \mod n$. \\

\item If $\G \cong A_5$ and $p\neq5$ then $G \cong A_{5}\times C_{n}$ or  $G$ is isomorphic to 
\begin{center}
$\left\langle \r, \sigma, \tau \right|\r^n=1,\sigma^2=\r^{\frac{n}{2}},\t^3=\r^{\frac{n}{2}},(\sigma\t)^5=\r^{\frac{n}{2}},\sigma\r\sigma^{-1}=\r,\t\r\t^{-1}=\r \rangle $
\end{center}

\item If $\G \cong U$ then $G \cong U \times C_n$ or  $G$ is isomorphic to 
\[ \<\r,\sigma_1,\dots ,\sigma_t|\r^n=\sigma_1^p= \dots =\sigma_t^p=1, \sigma_i\sigma_j=\sigma_j\sigma_i, \sigma_i\r\sigma_i^{-1}=\r^{l}, 1\leq i,j\leq t \>\]
where $(l,n)=1$ and $l^p \equiv 1 \mod n$.\\

\item If $\G \cong K_m$ then $G$ is isomorphic to  one of the following 

i) $
\< \r,\sigma_1,...,\sigma_t,v|\r^n=\sigma_1^p=...=\sigma_t^p=v^m=1, \sigma_i\sigma_j=\sigma_j\sigma_i,  v\r v^{-1}=\r, \sigma_i\r\sigma_i^{-1}=\r^{l}, \sigma_iv\sigma_i^{-1}=v^{k}, 1\leq i,j\leq t \> $ 
where $(l,n)=1$ and $l^p \equiv 1 \mod n, (k,m)=1$ and $k^p \equiv 1 \mod m$.  

ii) $ \< \r,\sigma_1, \dots ,\sigma_t \, | \, \r^{nm}  = \sigma_1^p = \dots = \sigma_t^p =1, \sigma_i\sigma_j = \sigma_j\sigma_i, \sigma_i \r \sigma_i^{-1} = \r^{l}, i \geq 1,  j\leq t \>$,   where  $(l, nm)=1$ and $l^p \equiv 1 \mod nm$. \\

\item If $\G \cong PSL_{2}(q)$ then $G\cong PSL_2(q)\times C_n$ or $SL_2(3)$. \\

\item If $\G \cong PGL_2 (q)$ then $G \cong PGL_2 (q) \times C_n$. \\
\end{enumerate}
\end{theorem}

\proof   See \cite{Sa1} for all the details. \\

For sake of completennes and also because of the fact that the signatures of $\Phi$   were crucial in determining all cases of the theorem above, we display all  these  signatures.  The proof can be found in \cite{Sa1}.  
\begin{lemma} \label{thm_3}
The signature of cover $\Phi : \X \to \X^G$ and dimension $\delta$ is given in Table~\ref{t2}, where $m=|PSL_2(q)|$ for cases 38-41 and $m=|PGL_2(q)|$ for cases 42-45.
\end{lemma}

%\small
\begin{table}[hb]
\begin{center}
\begin{tabular}{|c|c|c|c|}
\hline \hline
&&& \\
$Case$ & $\G$ & $\delta(G,C)$ & \textbf{C }$=(C_1,...,C_r)$ \\
%&&& \\
\hline \hline
$a$ &  & $\frac{g+n-1}{30(n-1)}-1$ & $(6,5,n,...,n)$  \\
$b$ &  & $\frac{g+5n-5}{30(n-1)}-1$ & $(6,5n,n,...,n)$  \\
$c$ & $A_5$  & $\frac{g+6n-6}{30(n-1)}-1$ & $(6n,5,n,...,n)$ \\
$d$ & & $\frac{g}{30(n-1)}-1$ & $(6n,5n,n,...,n)$  \\
\hline \hline
\end{tabular}
\vspace{0.4cm}
\caption{ $\delta$ for $\G \cong A_5$, $p=3$} \label{t3}
\end{center}
\end{table}

\normalsize

%\clearpage

\begin{remark} \label{r1}
The above Lemma gives signatures and dimensions for $p>5$. Since $\G \cong C_m, D_m, A_4, S_4, U, K_m, PSL(2,q), PGL(2,q)$ when $p=5$ and $\G \cong C_m, D_m, A_5$, $U$, $K_m$, $PSL(2,q)$, $PGL(2,q)$ when $p=3$, then all cases except $\G \cong A_5$ have ramification as $p>5$. However, $\G \cong A_5$ has different ramification. Hence, that case has signatures and dimensions as in Table~\ref{t3}.
\end{remark}

%\vspace{-0.5cm}

%\vspace{-1.2cm}

%\vspace{6ex}

\tiny

\begin{table}[ht]\label{tab2}
\begin{center}
\begin{tabular}{|c|c|c|c|c|}
\hline \hline 
&&&& \\
$\#$ & $\G$ & $\delta(G,C)$ & $\delta,n,g$ & \textbf{C} = $(C_1,...,C_r)$ \\
\hline \hline
&&&& \\
$1$ &$(p,m)=1$ & $\frac{2(g+n-1)}{m(n-1)}-1$ & $n<g+1$ & $(m,m,n,...,n)$  \\
$2$ & $C_m$ & $\frac{2g+n-1}{m(n-1)}-1$ & & $(m,mn,n,...,n)$ \\
$3$ &  & $\frac{2g}{m(n-1)}-1$ & $n<g$ &  $(mn,mn,n,...,n)$ \\
\hline \hline
$4$ & $(p,m)=1$ & $\frac{g+n-1}{m(n-1)}$ & & $(2,2,m,n,...,n)$  \\
$5$ &  & $\frac{2g+m+2n-nm-2}{2m(n-1)}$ & & $(2n,2,m,n,...,n)$  \\
$6$ & $D_{2m}$ & $\frac{g}{m(n-1)}$ & & $(2,2,mn,n,...,n)$ \\
$7$ & &$\frac{g+m+n-mn-1}{m(n-1)}$ & $n<g+1$ & $(2n,2n,m,n,...,n)$  \\
$8$ &  & $\frac{2g+m-mn}{2m(n-1)}$ & $g\neq 2$ & $(2n,2,mn,n,...,n)$  \\
$9$ &  & $\frac{g+m-mn}{m(n-1)}$ & $n<g$ & $(2n,2n,mn,n,...,n)$  \\
\hline \hline
$10$ &  & $\frac{n+g-1}{6(n-1)}$ & & $(2,3,3,n,...,n)$  \\
$11$ &$ A_4$ & $\frac{g-n+1}{6(n-1)}$ & & $(2,3n,3,n,...,n)$ \\
$12$ &  & $\frac{g-3n+3}{6(n-1)}$ & & $(2,3n,3n,n,...,n)$  \\
$13$ &  & $\frac{g-2n+2}{6(n-1)}$ & $\delta \neq 0$ & $(2n,3,3,n,...,n)$  \\
$14$ &  & $\frac{g-4n+4}{6(n-1)}$ & & $(2n,3n,3,n,...,n)$  \\
$15$ &  & $\frac{g-6n+6}{6(n-1)}$ & $\delta \neq 0$ & $(2n,3n,3n,n,...,n)$  \\
\hline \hline
$16$ &  & $\frac{g+n-1}{12(n-1)}$ & & $(2,3,4,n,...,n)$  \\
$17$ &  & $\frac{g-3n+3}{12(n-1)}$ & & $(2,3n,4,n,...,n)$  \\
$18$ &  &  $\frac{g-2n+2}{12(n-1)}$ & & $(2,3,4n,n,...,n)$ \\
$19$ &  & $\frac{g-6n+6}{12(n-1)}$ & & $(2,3n,4n,n,...,n)$  \\
$20$ &  $S_4$ & $\frac{g-5n+5}{12(n-1)}$ & & $(2n,3,4,n,...,n)$  \\
$21$ &  & $\frac{g-9n+9}{12(n-1)}$ & & $(2n,3n,4,n,...,n)$  \\
$22$ &  & $\frac{g-8n+8}{12(n-1)}$ & & $(2n,3,4n,n,...,n)$  \\
$23$ &  & $\frac{g-12n+12}{12(n-1)}$ & & $(2n,3n,4n,n,...,n)$  \\
\hline \hline
$24$ &  & $\frac{g+n-1}{30(n-1)}$ & & $(2,3,5,n,...,n)$  \\
$25$ &  & $\frac{g-5n+5}{30(n-1)}$ & & $(2,3,5n,n,...,n)$ \\
$26$ &  & $\frac{g-15n+15}{30(n-1)}$ & & $(2,3n,5n,n,...,n)$  \\
$27$ &  & $\frac{g-9n+9}{30(n-1)}$ & & $(2,3n,5,n,...,n)$  \\
$28$ &  $A_5$ & $\frac{g-14n+14}{30(n-1)}$ & & $(2n,3,5,n,...,n)$  \\
$29$ &  & $\frac{g-20n+20}{30(n-1)}$ & & $(2n,3,5n,n,...,n)$  \\
$30$ &  & $\frac{g-24n+24}{30(n-1)}$ & & $(2n,3n,5,n,...,n)$  \\
$31$ &  & $\frac{g-30n+30}{30(n-1)}$ & & $(2n,3n,5n,n,...,n)$  \\
\hline \hline
$32$ &  & $\frac{2g+2n-2}{p^t(n-1)}-2$ & & $(p^t,n,...,n)$  \\
$33$ & $U$ & $\frac{2g+np^{t}-p^t}{p^t(n-1)}-2$ & $(n,p)=1,n|p^t-1$ & $(np^t,n,...,n)$ \\
\hline \hline
$34$ &  & $\frac{2(g+n-1)}{mp^t(n-1)}-1$ & $(m,p)=1,m|p^t-1$ & $(mp^t,m,n,...,n)$  \\
$35$ &  & $\frac{2g+2n+p^t-np^t-2}{mp^t(n-1)}-1$ & $(m,p)=1,m|p^t-1$ & $(mp^t,nm,n,...,n)$  \\
$36$ & $K_m$ & $\frac{2g+np^t-p^{t}}{mp^t(n-1)}-1$ & $(nm,p)=1,nm|p^t-1$ & $(nmp^t,m,n,...,n)$\\
$37$ &  & $\frac{2g}{mp^t(n-1)}-1$ & $(nm,p)=1,nm|p^t-1$ & $(nmp^t,nm,n,...,n)$ \\
\hline \hline
$38$& & $\frac{2(g+n-1)}{m(n-1)}-1$ & $\left(\frac{q-1}{2},p\right)=1$ & $(\alpha,\beta,n,...,n)$  \\
$39$& $PSL_2(q)$  & $\frac{2g+q(q-1)-n(q+1)(q-2)-2}{m(n-1)}-1$ & $\left(\frac{q-1}{2},p\right)=1$ & $(\alpha,n\beta,n,...,n)$  \\
$40$& & $\frac{2g+nq(q-1)+q-q^2}{m(n-1)}-1$ & $\left(\frac{n(q-1)}{2},p\right)=1$ & $(n\alpha,\beta,n,...,n)$  \\
$41$& & $\frac{2g}{m(n-1)}-1$ & $\left(\frac{n(q-1)}{2},p\right)=1$ & $(n\alpha,n\beta,n,...,n)$  \\
\hline \hline
$42$& & $\frac{2(g+n-1)}{m(n-1)}-1$ & $(q-1,p)=1$ & $(2\alpha,2\beta,n,...,n)$  \\
$43$& $PGL_2(q)$  & $\frac{2g+q(q-1)-n(q+1)(q-2)-2}{m(n-1)}-1$ & $(q-1,p)=1$ & $(2\alpha,2n\beta,n,...,n)$  \\
$44$& & $\frac{2g+nq(q-1)+q-q^2}{m(n-1)}-1$ & $(n(p-1),p)=1$ & $(2n\alpha,2\beta,n,...,n)$ \\
$45$& & $\frac{2g}{m(n-1)}-1$ & $(n(q-1),p)=1$ & $(2n\alpha,2n\beta,n,...,n)$ \\ \\
\hline \hline
\end{tabular}
\medskip
\caption{The signature of curves and dimensions $\delta$ for $char >5 $} \label{t2}
\end{center}
\end{table}

\normalsize

\clearpage
%************************************************************************
\subsection{Equations of superelliptic curves}

Next we give the parametric equations of superelliptic curves based on their group of automorphisms.  Such equations for the first time were computed in \cite{Sa2}.   It is exactly the fact that their equations are easily determined that makes superelliptic curves quite attractive in applications.  Let   $\d$ be given as in Table~\ref{t2} and $M, \Lambda, Q, B, \Delta, \Theta$ and $\Omega$ are as follows: 

%\begin{small}
\begin{equation*}
\begin{split}
M = & \prod_{i=1}^\d \left(  x^{24}+\lambda_ix^{20}+(759-4\lambda_i)x^{16}+2(3\lambda_i+1228)x^{12} \right. \\
   & + \left. (759-4\lambda_i)x^8+\lambda_ix^4+1 \right) \\ \\
\Lambda= & \prod_{i=1}^\d \left(-x^{60}+(684-\lambda_i)x^{55}-(55\lambda_i+157434)x^{50}-(1205\lambda_i-12527460)x^{45}\right.\\
                &-(13090\lambda_i+77460495)x^{40}+(130689144-69585\lambda_i)x^{35}\\
                & +(33211924-134761\lambda_i)x^{30}+(69585\lambda_i-130689144)x^{25}\\
                & -(13090\lambda_i+77460495)x^{20}-(12527460-1205\lambda_i)x^{15}\\
                &  \left. -(157434+55\lambda_i)x^{10} +(\lambda_i-684)x^5-1 \right) \\ \\
Q  = & x^{30}+522x^{25}-10005x^{20}-10005x^{10}-522x^5+1\\ \\
B = &  \prod_{i=1}^\delta \displaystyle{\prod_{a \in H_t}}   \left( (x+a)-\lambda_i \right) \\ \\
\Theta = & \prod_{i=1}^\delta G_{\lambda_i} (x), \textit{ where  } G_{\lambda_i}= \left( x \cdot \prod_{j=1}^{\frac {p^t -1} m} (x^m-b_j) \right)^m  -\lambda_i \\ \\
\Delta = & \prod_{i=1}^{\d}(((x^q-x)^{q-1}+1)^{\frac{q+1}{2}}-\lambda_i(x^q-x)^{\frac{q(q-1)}{2}}) \\ \\
\Omega = & \prod_{i=1}^{\d}(((x^q-x)^{q-1}+1)^{q+1}-\lambda_i(x^q-x)^{q(q-1)}) \\ \\
\end{split}
\end{equation*}
%\end{small}

\noindent Then we have the following result. 

\begin{theorem}\label{thm_4}
Let $\X_g$ be e genus $g \geq 2$ algebraic curve defined over an algebraically closed field $k$, $G$ its automorphism group over $k$, and $C_n$  a cyclic normal subgroup of   $G$ such that $g (X_g^{C_n} ) =0$.  Then, the equation of $\X_g$ can be written as in one of the following cases as in Table~4.
\end{theorem}

%\vspace{-.5cm}
\small
\begin{table}[ht] \label{equations}
\begin{center}
\begin{tabular}{|c|c|c|}
\hline
%\hline
 & & \\
$\#$ &$ \bar G$& $y^n=f(x)$ \\
\hline
%\hline
1 &     & $x^{m\d}+a_1x^{m(\d-1)}+...+a_\d x^{m}+1$\\
2 & $C_m$    & $x^{m\d}+a_1x^{m(\d-1)}+...+a_\d x^{m}+1$\\
3 &     & $x(x^{m\d}+a_1x^{m(\d-1)}+...+a_\d x^{m}+1)$\\
\hline
%\hline
4 &        & $F(x):= \prod_{i=1}^\d(x^{2m}+\lambda_ix^m+1)$\\
5 &        & $(x^m-1)\cdot F(x)$\\
6 &        & $x\cdot F(x)$\\
7 &$D_{2m}$& $(x^{2m}-1)\cdot F(x)$\\
8 &        & $x(x^m-1)\cdot F(x)$\\
9 &        & $x(x^{2m}-1)\cdot F(x)$\\
\hline
%\hline
10 &     & $G(x):= \prod_{i=1}^\delta(x^{12}-\lambda_ix^{10}-33x^8+2\lambda_ix^6-33x^4-\lambda_ix^2+1)$\\
11 &     & $(x^4+2i\sqrt{3}x^2+1)\cdot G(x)$\\
12 &$A_4$& $(x^8+14x^4+1)\cdot G(x)$\\
13 &     & $x(x^4-1)\cdot G(x)$\\
14 &     & $x(x^4-1)(x^4+2i\sqrt{3}x^2+1)\cdot G(x)$\\
15 &     & $x(x^4-1)(x^8+14x^4+1)\cdot G(x)$\\
\hline
%\hline
16 &     & $M(x)$\\
17 &     & $  \left( x^8+14x^4+1 \right)  \cdot M(x)$\\
18 &     & $x(x^4-1) \cdot M(x)$\\
19 &     & $\left( x^8+14x^4+1 \right)  \cdot x(x^4-1) \cdot M(x)$\\
20 &$S_4$& $\left( x^{12}-33x^8-33x^4+1  \right)\cdot M(x)$\\
21 &     & $\left( x^{12}-33x^8-33x^4+1  \right)  \cdot \left( x^8+14x^4+1 \right)  \cdot M(x)$\\
22 &     & $\left( x^{12}-33x^8-33x^4+1  \right) \cdot x(x^4-1) \cdot M(x)$\\
23 &     & $\left( x^{12}-33x^8-33x^4+1  \right) \cdot \left( x^8+14x^4+1 \right)  \cdot x(x^4-1)  M(x)$\\
\hline
%\hline
24 &     & $\Lambda(x)$\\
25 &     & $x(x^{10}+11x^5-1) \cdot \Lambda(x)$\\
26 &     & $(x^{20}-228x^{15}+494x^{10}+228x^5+1)(x(x^{10}+11x^5-1))\cdot \Lambda(x)$\\
27 &     & $(x^{20}-228x^{15}+494x^{10}+228x^5+1)\cdot \Lambda(x)$\\
28 &$A_5$& $Q (x) \cdot \Lambda(x)$\\
29 &     & $x(x^{10}+11x^5-1).\psi(x)\cdot \Lambda(x)$\\
30 &     & $(x^{20}-228x^{15}+494x^{10}+228x^5+1)\cdot \psi(x)\cdot\Lambda(x)$\\
31 &     & $(x^{20}-228x^{15}+494x^{10}+228x^5+1)(x(x^{10}+11x^5-1))\cdot\psi(x)\cdot\Lambda(x)$\\
\hline
%\hline
32 &$U$& $B(x)$\\
33 &   & $B(x)$\\
\hline
34 &     & $\Theta(x)$\\
35 &$K_m$& $x\prod_{j=1}^{\frac{p^t-1}{m}}\left(x^m-b_j\right)\cdot\Theta(x)$\\
36 &     & $\Theta(x)$\\
37 &     & $x\prod_{j=1}^{\frac{p^t-1}{m}}\left(x^m-b_j\right)\cdot\Theta(x)$\\
\hline
%\hline
38 &          & $\Delta(x)$\\
39 &$PSL_2(q)$& $((x^q-x)^{q-1}+1)\cdot\Delta(x)$\\
40 &          & $(x^q-x)\cdot\Delta(x)$\\
41 &          & $(x^q-x)((x^q-x)^{q-1}+1)\cdot\Delta(x)$\\
\hline
%\hline
42 &          & $\Omega(x)$\\
43 &$PGL_2(q)$& $((x^q-x)^{q-1}+1)\cdot\Omega(x)$\\
44 &          & $(x^q-x)\cdot\Omega(x)$\\
45 &          & $(x^q-x)((x^q-x)^{q-1}+1)\cdot\Omega(x)$\\
\hline
%\hline
\end{tabular}
\vspace{.5cm}
\caption{The equations of the curves related to the cases in Table ~\ref{t2}}\label{t4}
\end{center}
\end{table}

\normalsize

%\vspace{-0.5cm}

\clearpage

Each case in the above table correspond to a $\delta$-dimensional family, where $\delta$ can be found in \cite{Sa2}. Moreover, our parametrizations are exact in the sense that the number of parameters in each case is equal to the dimension.  We would like to find invariants to classify isomorphism classes of these curves. 

%************************************************** 
 
\section{Isomorphism classes of superelliptic curves} 

A superelliptic curve $\X_g$ is given by an equation of the form $ y^n = f(x)$   for some degree $d$ polynomial $f(x)$.  Let us assume that 
\[ y^n = f(x)=\prod_{i=1}^s (x-\a_i)^{d_i}, \quad 0 < d_i < d.\]
We have that $\sum_{i=1}^s d_i =d$.  We call this the \textbf{standard form} of the curve.    The only places of $F_0=k(x)$ that ramify are the places which correspond to the points $x=\a_i$. We denote such places by $Q_1, \dots , Q_s$ and by $\B:=\{ Q_1, \dots , Q_s\}$ the set of these places.  The ramification indexes are $e (Q_i) = \frac n {(n , d_i)}$.

Hence, every set $\B$ determines a genus $g$ superelliptic curve $\X_g$. However, the correspondence between the sets $\B$ and the isomorphism classes of $\X_g$ is not a one to one correspondence. Obviously the set of roots of $f(x)$ does not determine uniquely the isomorphism class of $\X_g$ since every coordinate change in $x$ would change the set of these roots.  Such isomorphism classes are classified by the invariants of binary forms.  Invariants of binary forms of of degree up to eight are known by classical work of many invariant theorists and some more recent work, see \cite{vishi, Shi1, super4}. 
%However, for higher degree this is an open problem.  For some recent results see \cite{inv}. 

\subsection{Invariants of binary forms}

In this section we define the action of $ GL_2(k)$ on binary forms and discuss the basic notions of their invariants. Let $k[X,Z]$  be the  polynomial ring in  two variables and  let $V_d$ denote  the $(d+1)$-dimensional  subspace  of  $k[X,Z]$  consisting of homogeneous polynomials.
\begin{equation}  \label{eq1}
f(X,Z) = a_0X^d + a_1X^{d-1}Z + ... + a_dZ^d
\end{equation}
of  degree $d$. Elements  in $V_d$  are called  {\it binary  forms} of degree $d$.  We let $GL_2(k)$ act as a
group of automorphisms on $ k[X, Z] $   as follows:
\begin{equation}
 M =
\begin{pmatrix} a &b \\  c & d
\end{pmatrix}
\in GL_2(k), \textit{   then       }
\quad  M  \begin{pmatrix} X\\ Z \end{pmatrix} =
\begin{pmatrix} aX+bZ\\ cX+dZ \end{pmatrix}
\end{equation}
This action of $GL_2(k)$  leaves $V_d$ invariant and acts irreducibly on $V_d$.
\begin{remark}
It is well  known that $SL_2(k)$ leaves a bilinear  form (unique up to scalar multiples) on $V_d$ invariant. This
form is symmetric if $d$ is even and skew symmetric if $d$ is odd.
\end{remark}
Let $A_0$, $A_1$,  ... , $A_d$ be coordinate  functions on $V_d$. Then the coordinate  ring of $V_d$ can be
identified with $ k[A_0  , ... , A_d] $. For $I \in k[A_0, ... , A_d]$ and $M \in GL_2(k)$, define $I^M \in k[A_0,
... ,A_d]$ as follows
\begin{equation} \label{eq_I}
{I^M}(f):= I(M(f))
\end{equation}
for all $f \in V_d$. Then  $I^{MN} = (I^{M})^{N}$ and Eq.~(\ref{eq_I}) defines an action of $GL_2(k)$ on $k[A_0,
... ,A_d]$.
A homogeneous polynomial $I\in k[A_0, \dots , A_d, X, Z]$ is called a {\it covariant}  of index $s$ if
$$I^M(f)=\delta^s I(f),$$
where $\delta =\det(M)$.  The homogeneous degree in $a_1, \dots , a_n$ is called the {\it degree} of $I$,  and the
homogeneous degree in $X, Z$ is called the {\it  order} of $I$.  A covariant of order zero is called {\it
invariant}.  An invariant is a $SL_2(k)$-invariant on $V_d$.

We will use the symbolic method of classical theory to construct covariants of binary forms.    Let
$$f(X,Z):=\sum_{i=0}^n
\begin{pmatrix} n \\ i
\end{pmatrix}
a_i X^{n-i} \, Z^i, \quad  and \quad g(X,Z) :=\sum_{i=0}^m
  \begin{pmatrix} m \\ i
\end{pmatrix}
b_i X^{n-i} \, Z^i
$$
be binary forms of  degree $n$ and $m$ respectively with coefficients in $k$. We define the {\bf r-transvection}
$$(f,g)^r:= \frac {(m-r)! \, (n-r)!} {n! \, m!} \, \, \sum_{k=0}^r
(-1)^k
\begin{pmatrix} r \\ k
\end{pmatrix} \cdot
\frac {\partial^r f} {\partial X^{r-k} \, \,  \partial Z^k} \cdot \frac {\partial^r g} {\partial X^k  \, \,
\partial Z^{r-k} }
 $$
It is a homogeneous polynomial in $k[X, Z]$ and therefore a covariant of order $m+n-2r$ and degree 2. In general, the $r$-transvection of two covariants of order $m, n$ (resp., degree $p, q$) is a covariant of order $m+n-2r$  (resp., degree $p+q$).

For the rest of this paper $F(X,Z)$ denotes a binary form of order $d:=2g+2$ as below
\begin{equation}
F(X,Z) =   \sum_{i=0}^d  a_i X^i Z^{d-i} = \sum_{i=0}^d
\begin{pmatrix} n \\ i
\end{pmatrix}    b_i X^i Z^{n-i}
\end{equation}
where $b_i=\frac {(n-i)! \, \, i!} {n!} \cdot a_i$,  for $i=0, \dots , d$.  We denote invariants (resp., covariants) of binary forms by $I_s$ (resp., $J_s$) where the subscript $s$ denotes the degree (resp., the order).  We define the following covariants and invariants:
\begin{equation}
\begin{split}\label{covar}
\aligned
 I_2 & :=(F,F)^d,   \\
 I_4 & :=(J_4, J_4)^4,  \\
 I_6 & :=((F, J_4)^4, (F, J_4)^4)^{d-4},   \\
 I_6^\ast & :=((F, J_{12})^{12}, (F, J_{12})^{12})^{d-12},  \\
  M  & :=((F, J_4)^4, (F, J_8)^8)^{d-10}, \\
\endaligned
\qquad \aligned
& J_{4j}   :=   (F,F)^{d-2j}, \, \,  j=1, \dots , g, \\
& I_4'    :=     (J_8, J_8)^8, \\
& I_6^\prime  :=((F, J_8)^8, (F, J_8)^8)^{d-8}, \\
& I_3      :=(F, J_d)^d, \\
& I_{12}   :=(M, M)^8\\
\endaligned
\end{split}
\end{equation}
{\it Absolute invariants} are called $GL_2(k)$-invariants. We  define the following absolute invariants:
$$i_1:=\frac {I_4'} {I_2^2}, \, \,  i_2:=\frac {I_3^2} {I_2^3},
\, \,  i_3:=\frac {I_6^\ast  } {I_2^3}, \, \,  j_1 := \frac {I_6^{'}} {I_3^2}, \, \,   j_2:= \frac {I_6} {I_3^2},
\, \, u_1:=\frac {I_6^2} {I_{12}}, \, \, u_2:=\frac {(I_6^{'})^2} {I_{12}}
$$
$$
\v_1:= \frac {I_6} {I_6^\ast }, \, \,  \v_2:=\frac {(I_4^{'})^3} {I_3^4}, \, \, \v_3:= \frac {I_6} {I_6^{'}}, \,
\,  \v_4:=\frac {(I_6^\ast )^2} {I_4^3}.
$$
In the case $g=10$ and $I_{12}=0$ we define
\begin{equation}
\begin{split}
 I_6^\star & := ( (F, J_{16})^{16},  (F, J_{16})^{16} )^{d-16}), \\
S   \, \, & :=( J_{12}, J_{16} )^{12}, \\
I_{12}^\ast & := ( \, (J_{16}, S)^4, \, (J_{16}, S)^4 \, )^{12}\\
\end{split}
\end{equation}
and $\v_5\,  :=\frac {I_6^\star }{I_{12}^\ast}.$

For a given curve $\X_g$ we denote by $I(\X_g)$ or $i(\X_g)$ the corresponding  invariants.  When the above invariants are a good set of invariants to study the small genus curves, they are not a set of complete invariants for curves of arbitrary genus. 

\begin{exa} Let $C$ be a genus 4 curve with equation  
\[ y^3 = x^6 + a_5 x^5 + \dots + a_1 x + a_0,\]
defined over $\C$.  This curve has automorphism group $C_3$.  The family $V$ of such curves is a 3-dimensional variety.  The isomorphism classes of curves in this variety are determined by Igusa invariants $J_2, J_4, J_6, J_{10}$, see \cite{vishi, sh_03} for their definitions.  Two curves $C$ and $C^\prime$ in V are isomorphic if and only if 
\[ \left(J_2 (C), J_4 (C), J_6 (C), J_{10} (C)\right) = \lambda \cdot (J_2 (C^\prime) , J_4 \left(C^\prime) , J_6 (C^\prime), J_{10} (C^\prime) \right) \]
for some $\lambda \neq 0$.
\end{exa}

\begin{lemma}
Let $\X_g$ be a superelliptic curves of genus $g \geq 2$.  The following statements are true.   

i) If $\G \iso A_4$ then $I_4 (\X_g) =0$.   

ii)  If $\G \iso A_5$ then $(J_i, J_i )^i =0 $ for $i=4, 8, 16, 28$. 
\end{lemma}

\proof 
See \cite{super4}  for the proof of these and other properties of superelliptic curves in terms of invariants of binary forms. 

\endproof

%************************************************************************************
\section{\ss-invariants of superelliptic curves}  

In this section we will introduce \ss-invariants of superelliptic curves.  These invariants were introduced in \cite{g_sh} for hyperelliptic curves and generalized in  \cite{AK} for superelliptic curves.  Here we simply follow the approach from \cite{AK}. 

Let $k$ be an algebraic closed field of characteristic $p \geq 0$. Let $F_0 = k(x)$ be the function field of the projective line $\P$. We consider a cyclic extension of $F_0$ of degree $n$ of the form $F:=k(x,y)$ where 
\begin{equation}
y^n= \prod_{i=1}^{s}(x-\rho_i)^d_{i} =: f(x), \quad o<d_i<n.
\end{equation}

If $d:= \sum_{i=1}^{s}d_i \equiv  0\mod n$ then the place at infinity does not ramify at the above extension. The only places at $F_0$ that are ramified are the places $P_i$ that correspond to the points $x=\rho_i$ and the corresponding ramification indices are given by 
\[ e_i = \frac{n}{(n,d_i)}. \]
Moreover if $(n,d_i) =1$ then the places $P_i$ are ramified completely and the Riemann-Hurwitz  formula implies that the function field $F$ has genus
\[ g= \frac{(n-1)(s-2)}2. \]
Notice that the condition $g\geq 2$ is equivalent to $s\geq 2\frac{n+1}{n-1}$. In particular, $s >2$.  

For the proof of the following Lemmas se \cite{AK}. 

\begin{lemma}
Let $G= Aut(F)$. Suppose that a cyclic extension $F/F_0$ of the rational function field $F_0$ is ramified completely at $s$ places and $n:= |Gal(F/F_0)|$. If $2n <s$ then $Gal(F/F_0)\normal G$.
\end{lemma}

\begin{lemma}
Suppose that $\tau$ is an extra automorphism of $F$, and let $s$ be the number of ramified places at the extension $F/F_0$ and let $d$ be the degree of the defining polynomial. Then $\delta|s, \delta|d$ and the defining equation of $F$ can be written as
\[y^n = \sum_{i-0}^{d/\delta}a_ix^{\delta \cdot i},    \]
where $a_0 =1$. 
\end{lemma}

We will say that the superelliptic curve is  in \textbf{normal form} if and only if it is given by an equation:
\[y^n = x^s + \sum_{i-1}^{\frac{d}{\delta}}a_ix^{\delta \cdot i} + 1.    \]
%

%%%%%%%%%%%Valmira

Parametrizing superelliptic curves  that admit an extra automorphism of order $\d$, is the set of coefficients $\{ a_{s / \d-1}, \cdots, a_1\}$ of  a normal form up to a change of coordinate in $x$.   The condition $\tau (x)= \zeta x$, implies that $\bar{\tau}$ fixes the places $0, \infty$. Moreover we can change the defining equation by a morphism $\g \in PGL(2,k)$ of the form $ \g: x \rightarrow mx$ or $ \g: x \rightarrow \frac{m}{x}$  so that the new equation is again in normal form.  Substituting  
$a_0= (-1)^{d/s} \prod_{i=1}^{d/s} \b_i^s$
we have
\[(-1)^{s/ \d} \prod_{i=1}^{s/ \d} \g(\b_i)^\d =1\]
and this gives $m^s= (-1)^{s/\d}$. Then,  $x$ is determined up to a coordinate change by the subgroup $D_{s/ \d}$ generated by
\[ \t_1: x \rightarrow \epsilon  x, \thinspace \t_2 : x \rightarrow \frac{1}{x} \]
where $\epsilon$ is a primitive $s/\d$-root of one, see \cite{g_sh} for details. 

The action of $D_{s/\d}$ on the parameter space $k(a_1,\dots, a_{s/\d} )$ is given by 
\[
\begin{split}
 \t_1 : &  \, a_i  \rightarrow \epsilon^{\d i}  a_i,   \textit{ for } i= 1, \dots s/\d \\
 \t_2 : & \, a_i  \rightarrow   a_{d/ \d-i},       \textit{ for } i= 1,  \dots [s/\d]  \\
 \end{split}
 \]

Notice that if $s/\d=1$ then the above actions are trivial, therefore the normal form determines the equivalence class.  If $s/\d=2$ then 
$$ \t_1(a_1) = -a_1,   \tau_1 (a_2)= a_2, \tau_2= 1 $$
and the action is not dihedral but cyclic on the first vector.

\begin{lemma}  Assume that $s/\d> 2$. The fixed field $k(a_1,a_2,\cdots a_{s/\d})^{D_{s/\d}}$ is the same as the function field of the variety $\L_{n,s,\d}$.
\end{lemma}

\proof See \cite{AK} for the proof. 

\begin{lemma} Let $r:= s/\d >2$ The elements 
$$ \u_i:= a_1^{r-i}a_1 + a_{r-1}^{r-i}a_{r-i},  \textit{ for } i=1, \dots ,r$$
 are  invariants under the action of the group $D_{s/\d}$  defined as above. 
\end{lemma}

\proof See \cite{AK} for the proof. 

The elements $\u_i$ are called the \textbf{dihedral invariants} or \textbf{\ss-invariants}  of  $D_{s/\d}$.

 \begin{theorem}
 Let $\u=(\u_1, \dots , \u_r )$ be the $r-$ tuple of  \ss-invariants. Then $k(\L_{s,n,\d})= k(\u_1, \dots, \u_r)$.
 \end{theorem}

%%%%%%%%%%%%%%%%%%%%%%%%%%%%%%%%%%%%%

\begin{exa} For genus $g=2$ all curves are hyperelliptic and therefore superelliptic.  The generic curve of genus 2 is given by $y^2 = f(x)$, where $\deg (f) = $ 5 or 6. The space is determined by the invariants of binary sextics. When the curves have extra automorphisms then we have two main cases. \\

i) The first case is when there is an automorphism of degree 5. Then the full automorphism group $G$ is isomorphic to $\Z_{10}$ and the corresponding space has dimension zero.  There is only one curve in this case (up to isomorphism), which is given by $y^2=x^5-1$. \\

ii) In the second case, the extra automorphism is an involution. Then, $G$ is isomorphic to the Klein four-group $V_4$ and the curve has equation
\[ y^2 = x^6 + a_1 x^4 + a_2 x^2 + 1.\]
The  \ss-invariants are
\begin{equation}
\s_1=a_1^3 +a_2^3, \quad   \quad  \s_2=2 a_1 a_2, 
\end{equation}
see \cite{sh_2000} for a detailed study of this case.   In \cite{sh_2000} was the first time that such invariants were defined and later generalized in \cite{g_sh}. 
\end{exa}

%****************************************************
\section{Computational aspects of invariants of superelliptic curves}

In this part we give a quick introduction to some computational aspects of \ss-invariants.  A more detailed study of superelliptic curves and their computational invariants will appear in \cite{team2}. \\

\noindent \textbf{Problem}  Given a genus $g\geq 3$.   

1) Find the lattice of inclusions of all the cases based on the automorphism groups. 

2) Compute relations among \ss-invariants for every group of the table. \\

In other words, we would like to characterize for every group $G$ the locus of the curves in each case in the Table 2, in terms of invariants of these curves and determine the inclusions among such loci. 
While such lattice can be computed using only   group theory methods, from the computational viewpoint this is really not very useful. Instead such lattice and such loci need to be computed in terms of coefficients of the curves, or more precisely invariants of the curves.  A step further would be to characterize the Jacobians of curves in these loci. This can be done through the theory of theta functions as in 
\cite{beshaj_theta}. 

%\subsection{Computing \ss-invariants and invariants of binary forms.}

\subsection{A Maple package for computing with superelliptic curves}

Computing the \ss-invariants we first need the equation of the curve in the normal form \[ y^n = f(x). \]   
Once the normal form is determined then it is rather straight forward to compute the \ss-invariants.  We have implemented some of these tasks in Maple and display the codes below.   \\

\begin{verbatim}

normalpol:=proc(f,x)      # Computes the normal form of a polynomial.

    local a,n,f1;    
    n:=degree(f,x); f1:=f/coeff(f,x,0);  a:=coeff(f1,x,n)^(1/n);
    
RETURN(subs(x=x/a,f1));
end:

s_inv := proc(f, x)            # Computing the s-invariants. 

       local i, a, g, s;      
       
     g:=(degree(f, x) - 2)/2;

     for i to g do 
        a[i]:=coeff(f, x, 2*i) 
     od:    
     
     for i to g do 
         s[i]:=factor(a[1]^(g-i+1)*a[i]+a[g]^(g-i+1)*a[g-i+1]) 
     od;
     
RETURN([seq(s[i],i=1..g)]);
end;

fg_s:=proc (f,g,x,y,s)       # Computing the s-transvection of 
    local n,m,fg,k;              # binary forms f and g. 

  n:=degree(f,{x,y});
  m:=degree(g,{x,y});
  
  fg:=(n-s)!*(m-s)!/(n!*m!)*add((-1)^k*( s!/(k!*(s-k)!))
       *diff(f,x$(s-k),y$k) *diff(g,x$k,y$(s-k)),k=0..s );
       
  RETURN(expand(fg));
end:

fg_s2:=proc (f,g,x,y)    # fg_s2(f,g,x,y) := fg_s(f,g,x,y,deg(f)) 
   local n,f2,g2,k;         # deg(g)=deg(f) 

  n:=degree(f,{x,y});    # returns an invariant (which means order=0)
  
  f2:=collect(f,[x,y]);
  g2:=collect(g,[x,y]);
  
  1/(n!*n!)*add((-1)^k*( n!/(k!*(n-k)!))*(n-k)!*k!
   *coeff(coeff(f2,x,n-k),y,k)
   * (n-k)!*k!*coeff(coeff(g2,x,k),y,n-k) ,k=0..n );
   
  RETURN(expand(%));
end:

homogpol:=proc(f,x,y)    # Converts a polynomial to a homogenous one.
  RETURN(expand(subs(x=x/y,f*y^degree(f,x))));
end:

J_i:=proc(F,x,y,i)
  fg_s(F,F,x,y,degree(F,{x,y})-i/2);
  
  RETURN(%);
end:

I4prime:=proc(F,x,y)
  J_i(F,x,y,8);
  fg_s2(%,%,x,y);
  
  RETURN(%);
end:

I2:=proc(F,x,y)
  fg_s2(F,F,x,y);
  
  RETURN(%);
end:

I3:=proc(F,x,y)
  J_i(F,x,y,degree(F,{x,y}));
  fg_s2(F,%,x,y);
  
  RETURN(%);
end:
\end{verbatim}

%********************************************************************

%\clearpage
\section{Genus 3} 

In this section we will determine all the superelliptic curves of genus 3. Completing the case in positive characteristic is a natural extension of the methods used here.  

%The case for $g=3$ has appeared before in the literature for the case of characteristic zero, see \cite{GSS}.  

\subsection{Automorphism groups of genus 3 superelliptic curves}

Applying Thm.~\ref{thm_2}  we obtain the automorphism groups of a genus 3 superelliptic curves  defined over algebraically closed field of characteristic $p \neq 2$. Below we list GAP group ID's of those groups. 

%\begin{figure}[hbt!]
%\includegraphics[width=8.0cm]{locichart}
%\caption{Inclusions among the loci for genus 3 curves}
%\end{figure}
%

\begin{small}
\begin{center}
\begin{figure}[ht!]
\begin{tikzpicture} [scale=.65] %[thick]
\node at (-1, 0) [rectangle,draw=black!50,  fill=yellow!20] { 0 };
\node at (-1, 2) [rectangle,draw=black!50,  fill=yellow!20] { 1 };
\node at (-1, 4) [rectangle,draw=black!50,  fill=yellow!20] { 2 };
\node at (-1, 6) [rectangle,draw=black!50,  fill=yellow!20] { 3 };
\node at (-1, 8) [rectangle,draw=black!50,  fill=yellow!20] { 4 };
\node at (-1, 10) [rectangle,draw=black!50,  fill=yellow!20] { 5 };
%\node at (0, 13) [rectangle,draw=black!50,  fill=yellow!20] { dimension $\delta$ };
%%%
\node (node1)  at (1, 10) [ draw=black!50,  fill=blue!60] { 1, $C_2$ };

\node (node12)  at (10, 8) [draw=black!50,  fill=green!60] { 12, $C_2$ };

\node (node13)  at (12, 6) [draw=black!50,  fill=green!60] { 13, $V_4$ };
\node (node2)  at (2, 6) [draw=black!50,  fill=blue!60] { 2, $V_4$ };

\node (node3)  at (3.5, 4) [draw=black!50,  fill=blue!60] { 3, $C_2^3$ };
\node (node4)  at (5.5, 4) [draw=black!50,  fill=blue!60] { 4, $C_4$ };
\node (node14)  at (8, 4) [draw=black!50,  fill=yellow!60] { 14, $C_3$ };
\node (node15)  at (10, 4) [draw=black!50,  fill=green!60] { 15, $S_3$ };
\node (node16)  at (13, 4) [draw=black!50,  fill=green!60] { 16, $D_8$ };

\node (node5)  at (1, 2) [draw=black!50,  fill=blue!60] { 5 };
\node (node6)  at (2.5, 2) [draw=black!50,  fill=blue!60] { 6 };
\node (node7)  at (4.5, 2) [draw=black!50,  fill=blue!60] { 7 };
\node (node17)  at (8, 2) [draw=black!50,  fill=yellow!60] { 17, $C_6$ };
\node (node18)  at (10, 2) [draw=black!50,  fill=yellow!60] { 18 };
\node (node19)  at (13, 2) [draw=black!50,  fill=green!60] { 19, $S_4$ };

\node (node8)  at (0.5, 0) [draw=black!50,  fill=blue!60] { 8, $C_{14}$ };
\node (node9)  at (2, 0) [draw=black!50,  fill=blue!60] { 9 };
\node (node10)  at (3, 0) [draw=black!50,  fill=blue!60] { 10 };

\node (node11)  at (4.5, 0) [draw=black!50,  fill=blue!60] { 11 };

\node (node20)  at (6, 0) [draw=black!50,  fill=yellow!60] { 20, $C_9$ };
\node (node21)  at (8, 0) [draw=black!50,  fill=yellow!60] { 21 };
\node (node22)  at (11, 0) [draw=black!50,  fill=yellow!60] { 22 };
\node (node23)  at (13, 0) [draw=black!50,  fill=green!60] { 23, $L_3(2)$ };
 
%%%%%%%%%%   drawing nodes
\draw [->, blue!80] (node1.south) -- (node2.north);
\draw [->, blue!80] (node1.south) -- (node4.north);
\draw [->, blue!80] (node1.south) -- (node8.north);

\draw [->, blue!80] (node2.south) -- (node3.north);
\draw [->, blue!80] (node2.south) -- (node5.north);
\draw [->, blue!80] (node2.south) -- (node6.north);

\draw [->, blue!80] (node3.south) -- (node7.north);
\draw [->, blue!80] (node4.south) -- (node7.north);

\draw [->, blue!80] (node5.south) -- (node9.north);

\draw [->, blue!80] (node6.south) -- (node9.north);
\draw [->, blue!80] (node6.south) -- (node11.north);

\draw [->, blue!80] (node7.south) -- (node10.north);
\draw [->, blue!80] (node7.south) -- (node11.north);

\draw [->, green!80] (node12.south) -- (node2.north);
\draw [->, green!80] (node12.south) -- (node17.north);
\draw [->, green!80] (node12.south) -- (node15.north);
\draw [->, green!80] (node12.south) -- (node13.north);

\draw [->, green!80] (node13.south) -- (node16.north);

\draw [->, green!80] (node14.south) -- (node6.north);
\draw [->, green!80] (node14.south) -- (node20.north);
\draw [->, green!80] (node14.south) -- (node17.north);
\draw [->, green!80] (node15.south) -- (node6.north);
\draw [->, green!80] (node15.south) -- (node19.north);
\draw [->, green!80] (node16.south) -- (node7.north);
\draw [->, green!80] (node16.south) -- (node18.north);
\draw [->, green!80] (node16.south) -- (node19.north);
\draw [->, green!80] (node17.south) -- (node21.north);
\draw [->, green!80] (node18.south) -- (node21.north);

\draw [->, green!80] (node19.south) -- (node11.north);
\draw [->, green!80] (node19.south) -- (node22.north);
\draw [->, green!80] (node19.south) -- (node23.north);
\end{tikzpicture}
\caption{The lattice of genus 3 case. The blue items correspond to hyperelliptic curves, the yellow ones to the other superelliptic cases. }
\end{figure}
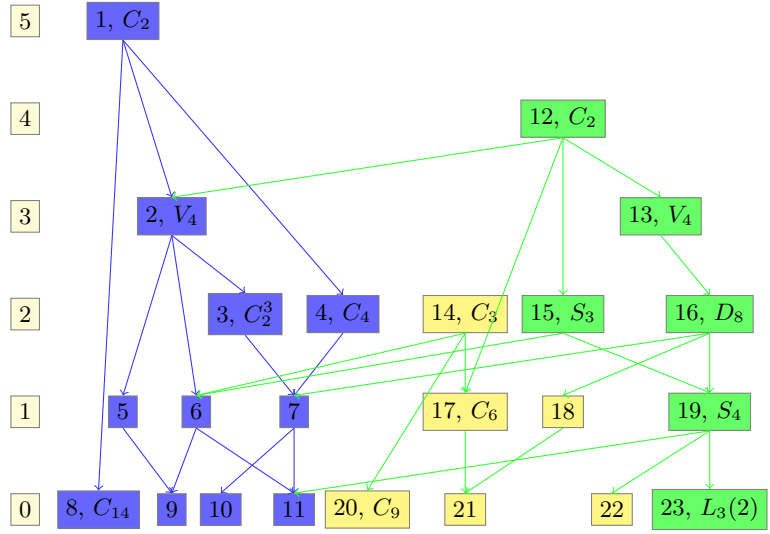
\end{center}
\end{small}

\begin{lemma}\label{thm_g_3}
Let $\X_g$ be a genus 3 superelliptic curve  defined over a field of characteristic $p \neq 2$. Then the automorphism groups of $\X_g$ are as follows.

\begin{description}
%    \item[i)] $p=0$ : $(2,1)$, $(4,2)$, $(3,1)$, $(4,1)$, $(8,2)$, $(8,3)$, $(7,1)$, $(21,1)$, $(14,2)$, $(6,2)$, $(12,2)$,
 %   $(9,1)$, $(8,1)$, $(8,5)$, $(16,11)$, $(16,10)$, $(32,9)$, $(30,2)$, $(42,3)$, $(12,4)$, $(16,7)$, $(24,5)$, $(18,3)$,
%    $(16,8)$, $(48,33)$, $(48,48)$.

\item[i)] $p=3$: $(2,1)$, $(4,2)$, $(3,1)$, $(4,1)$, $(8,2)$, $(8,3)$, $(7,1)$, $(14,2)$, $(6,2)$, $(8,1)$, $(8,5)$, $(16,11)$, $(16,10)$, $(32,9)$, $(30,2)$, $(16,7)$, $(16,8)$, $(6,2)$.

\item[ii)] $p=5$: $(2,1)$, $(4,2)$, $(3,1)$, $(4,1)$, $(8,2)$, $(8,3)$, $(7,1)$, $(21,1)$, $(14,2)$, $(6,2)$, $(12,2)$,
    $(9,1)$, $(8,1)$, $(8,5)$, $(16,11)$, $(16,10)$, $(32,9)$, $(42,3)$, $(12,4)$, $(16,7)$, $(24,5)$, $(18,3)$,
    $(16,8)$, $(48,33)$, $(48,48)$.

    \item[iii)] $p=7$: $(2,1)$, $(4,2)$, $(3,1)$, $(4,1)$, $(8,2)$, $(8,3)$, $(7,1)$, $(21,1)$, $(6,2)$, $(12,2)$,
    $(9,1)$, $(8,1)$, $(8,5)$, $(16,11)$, $(16,10)$, $(32,9)$, $(30,2)$, $(42,3)$, $(12,4)$, $(16,7)$, $(24,5)$, $(18,3)$,
    $(16,8)$, $(48,33)$, $(48,48)$.

\item[iv)] $p=0$ or $p > 7$: $(2,1)$, $(4,2)$, $(3,1)$, $(4,1)$, $(8,2)$, $(14,2)$, $(6,2)$,  
    $(9,1)$, $(8,5)$, $(16,11)$,  $(32,9)$,   $(12,4)$, $(16,13)$, $(24,5)$, 
      $(48,33)$, $(48,48)$, (96, 64).
\end{description}
\end{lemma}

Recall that the list for  $p=0$ is the same as for $p > 7$.   In the diagram below we display the inclusion among the loci in the case of genus 3.    We will briefly discuss the superelliptic curves and display their equations. 

%*********************************************************************
\subsection{Equations of hyperelliptic curves of genus three}

Let $\X_3$ be a hyperelliptic curve of genus 3.  In Tab.~\ref{t5} we list the automorphism groups of genus 3 hyperelliptic curves. The first column is the case number, in the second column the groups which occur as full automorphism groups are given, and the third
column indicates the reduced automorphism group for each case. The dimension $\delta$ of the locus and the
equation of the curve are also  are given in the next two columns. The last column is the {\sc GAP}
identity of each group in the library of small groups in {\sc GAP}.

\begin{table*} [hbt!]
\begin{small}
 \vskip 0.2cm
\begin{center}
\begin{tabular}{||c|c|c|c|c|c||}
\hline \hline & &   & & &\\
 &$\Aut(\X_g)$ & $\bAut(\X_g)$ & $\, \delta \, $ &equation  $y^2= f(x) $  &  Id.\\
%&         &    & &   &\\
%
\hline \hline  &&&& &  \\
 1 & $\Z_2$ &$\{1\}$ & 5&$x(x-1)(x^5+ax^4+bx^3+cx^2+dx+e)$&  $(2,1)$ \\
%
%\hline
& & & &   &\\
2 & $\Z_2 \times \Z_2$ &$\Z_2$& 3&$ x^8 + a_3 x^6 + a_2 x^4 + a_1 x^2 + 1$ &  $(4,2)$\\
3 & $\Z_2^3$ &$D_4$ &2&$(x^4+ax^2+1)(x^4+bx^2+1)$   &$(8,5)$\\
4 & $\Z_4$ & $\Z_2$ &2&$x(x^2-1)(x^4+ax^2+b)$&  $(4,1)$\\
5 & $\Z_2\times \Z_4$ &$D_4$ & 1&$(x^4-1)(x^4+ax^2+1)$&    $(8,2)$\\
6 & $D_{12}$ &$D_6$  &1&$x(x^6+ax^3+1)$   &$(12,4)$\\
7 & $\Z_2 \times D_8$ &$D_8$ &1&$x^8+ax^4+1$ &    $(16,11)$ \\
8 & $\Z_{14}$ &$\Z_7$ &0&$x^7-1$ &  $(14,2)$ \\
%
%\hline
& & & &  &\\
9 & $U_6$ &$D_{12}$ & 0&$x(x^6-1)$    & $(24,5)$\\
10 & $V_8$ &$D_{16}$  &0& $x^8-1$&   $(32,9)$\\
%
%
%\hline
& & & &   &\\
11 & $\Z_2 \times S_4$ & $S_4$  &0  & $x^8+14x^2+1$ &   $(48,48)$ \\
&  & & &   &\\
 \hline\hline
\end{tabular}
\end{center}
\end{small}
\caption{$\Aut(\X_3)$ for hyperelliptic $\X_3$}\label{t5}
\end{table*}

\medskip

%  Case 1. 

\noindent{\textbf{Case 1:  $\Z_2$-hyperelliptic: }   Then, the equation of $\X_3$ is given by \[ y^2 = f(x)\] where $\deg f =$ 7 or 8.  
To have an explicit way of describing a point in the moduli space of hyperelliptic curves of genus 3 we need
the generators of the field of invariants of binary octavics. These invariants are described in terms of
covariants of binary octavics. Such covariants were first constructed by van Gall who showed that  the graded
ring of covariants is generated by  70 covariants and explicitly constructed them,  see \cite{Shi1}.   

Let  $f(X,Y)$ be the binary octavic
\[f(X,Y) = \sum_{i=0}^8 a_i X^i Y^{8-i}.\]
We define the following covariants:
\begin{equation*}
\begin{split}
 & g=(f,f)^4, \quad k=(f, f )^6, \quad h=(k,k)^2, \quad m=(f,k)^4, \\
 & n=(f,h)^4, \quad p=(g,k)^4, \quad q=(g, h)^4. \\
\end{split}
\end{equation*}

Then the following
\begin{equation}\label{J}
\begin{split}
& J_2=(f,f)^8,\quad J_3=(f,g )^8,\quad J_4=(k,k)^4,  \quad   J_5=(m,k)^4,\\
&  J_6 = (k,h )^4, \quad  J_7=  (m,h)^4  \quad J_8=(p,h)^4, \quad J_9=(n,h)^4, \quad J_{10}=(q,h)^4    \\
\end{split}
\end{equation}
are $SL_2(k)$-invariants. Shioda has shown that the ring of invariants is a finitely generated module of
$k[J_2, \dots , J_7]$, see \cite{Shi1} for more details.  \\

%%%%%%%  Case V_4 - hyperelliptic 

\noindent{\textbf{Case 2:  $V_4$-hyperelliptic: }     Then,  $\X_3$ has normal equation
\[Y^2= X^8 + a_3 \, X^6 + a_2 \, X^4 + a_1 \, X^2 + 1,\]
see  \cite{g_sh}. The \ss-invariants of $\X_3$ are
\[\s_1   = a_1^4+a_3^4, \quad \s_2   =(a_1^2+a_3^2)\, a_2, \quad \s_3   =\, 2 \, a_1 \, a_3.\]
If $a_1=a_3=0$,   then $\s_1=\s_2=\s_3=0$. In this case  \[w := a_2^2\] is invariant. Thus, we define
\begin{equation} \label{def_u}
\begin{split}
\u (\X_3) \, = \left\{ \aligned & w & \mathrm{if} \quad a_1=a_3=0,\\
 & (\s_1, w, \s_3) & \mathrm{if} \quad a_1^2+a_3^2=0\  \mathrm{and} \, \, a_2\neq 0,\\
 & (\s_1, \s_2, \s_3) & \mathrm{otherwise}. \\ \endaligned
\right.
\end{split}
\end{equation}

The expressions of these covariants are very large in terms of the
coefficients of the curve and difficult to compute. However, in terms of the \ss-invariants $\s_1, \s_2,
\s_3$ these expressions are smaller. Analogously, $J_{14}$ is the discriminant of the octavic.  All these invariants are determined explicitly in terms of the \ss-invariants in \cite{GSS}.

We denote by  $\L_3$ the sublocus of $\M_3$ of hyperelliptic curves with automorphism group $V_4$.  This is a closed subvariety of $\M_3$  determined as below as shown in \cite{GSS}. The following are true, see \cite{GSS} for their proofs.

\begin{remark}

i)  $k(\L_3) = k(\s_1, \s_2 , \s_3)$. \label{thm_g_4}
 
 ii)  The relations among the \ss-invariants for other hyperelliptic curves of genus 3 are given in the Fig.~\ref{s-invariants_g_3}.
 \end{remark} 
 
 \bigskip

\begin{small}
\begin{center}
\begin{figure}[hbt!]
\begin{tikzpicture} [scale=.65] %[thick]
\node at (-2, 0) [rectangle,draw=black!50,  fill=yellow!20] { 0 };
\node at (-2, 2) [rectangle,draw=black!50,  fill=yellow!20] { 1 };
\node at (-2, 4) [rectangle,draw=black!50,  fill=yellow!20] { 2 };
\node at (-2, 6) [rectangle,draw=black!50,  fill=yellow!20] { 3 };

%\node at (0, 13) [rectangle,draw=black!50,  fill=yellow!20] { dimension $\delta$ };
%%%

\node (node2)  at (8, 6) [draw=black!50,  fill=blue!60] { 2, $V_4$: $(\s_1, \s_2, \s_3)$ };

\node (node3)  at (14, 4) [draw=black!50,  fill=blue!60] { 3, $C_2^3$, $2\s_1=\s_3^2$ };

\node (node5)  at (2, 2) [draw=black!50,  fill=blue!60] { 5, $\Z_2 \times \Z_4$, $2\s_1=-\s_3^2$ };
\node (node6)  at (8, 2) [draw=black!50,  fill=blue!60] { 6, $D_{12}$, Eq.~\eqref{d_u_12} };
\node (node7)  at (14, 2) [draw=black!50,  fill=blue!60] { 7, $\Z_2 \times D_8$, $a_1=a_3$ };

%\node (node8)  at (2, 0) [draw=black!50,  fill=blue!60] { 8, $C_{14}$ };
\node (node9)  at (2, 0) [draw=black!50,  fill=blue!60] { 9, $\s_2=0$ };

\node (node11)  at (8, 0) [draw=black!50,  fill=blue!60] {11,  $(0, 196, 0) \textrm{\ or\ }\left(\frac {8192} {81}, - \frac {1280} {27}, \frac {128} 9\right)$ };

\node (node10)  at (15, 0) [draw=black!50,  fill=blue!60] {10, $a_1=a_2=a_3=0$ };
 
%%%%%%%%%%   drawing nodes

\draw [->, blue!80] (node2.south) -- (node3.north);
\draw [->, blue!80] (node2.south) -- (node5.north);
\draw [->, blue!80] (node2.south) -- (node6.north);

\draw [->, blue!80] (node3.south) -- (node7.north);

\draw [->, blue!80] (node5.south) -- (node9.north);

\draw [->, blue!80] (node6.south) -- (node9.north);
\draw [->, blue!80] (node6.south) -- (node11.north);

\draw [->, blue!80] (node7.south) -- (node10.north);
\draw [->, blue!80] (node7.south) -- (node11.north);
\end{tikzpicture}
\caption{Relations among \ss-invariants for hyperelliptic curves of genus 3 with extra involutions.}
\label{s-invariants_g_3}
\end{figure}
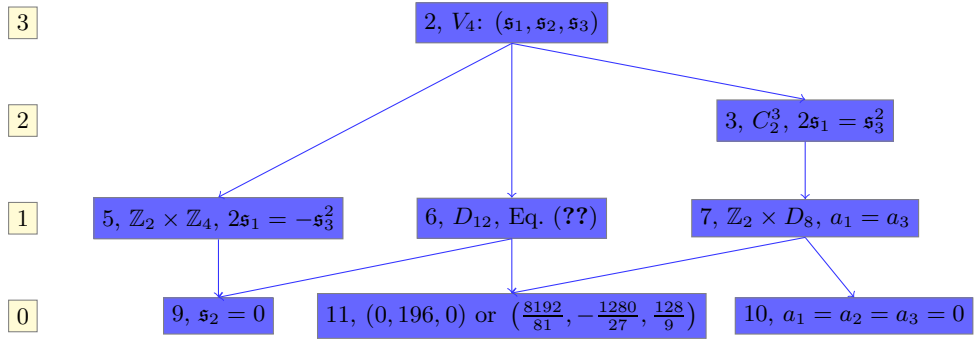
\end{center}
\end{small}

\bigskip  
%\clearpage
%***************************************************************************
\subsection{Equations of other superelliptic curves of genus 3.}  

In this section we take a quick glance of  all superelliptic curves of genus 3 in all positive characteristics $\neq 2$.  Similar tables are computed for all genus $g \leq 10$ in all characteristics and include for each curve the normal equation of the curve, the automorphism group, the invariants of the corresponding binary forms, $\s$-invariants, the dimension of the corresponding moduli space. In a current project the half-integer theta characteristics will be computed in each case and the equation of the corresponding curve in terms of these characteristics.

In the Table below we present these curves for $p=7$ for the cases 1-32 of the Table~2 so we can give an idea how this tables will look on the website with all the data.  Clicking on each curve will display all the information about the curve such as the automorphism group, invariants of the ninary form, \ss-invariants, an equation of the curve in terms of the theta-nulls, etc. 
      
        \begin{table}[] \label{genus 3_char7}
      \caption{Superelliptic curves of genus three}
      \vskip 0.5cm 
      \vskip 0.5cm
      \begin{center}
      \renewcommand{\arraystretch}{1.26}
      \begin{tabular}{||c|c|c|c|c||}
      \hline
      \hline
      $\#$ & Case & $\G$ & $n$ & Equation  \\
      \hline
      \hline      
        \multicolumn{5}{c}{Genus 3, $p = 7$} \\
      \hline
      \hline
             {\bf 1} &  1 &    &   &  $x^7+ a_1x^6 + a_2x^5 + a_3x^4 +a_4x^3 +a_5x^2 +a_6x +1$  \\
                {\bf 2} & 4 &    &   &  $(x^2 +a_1x +1)(x^2+a_2x +1)(x^2 +a_3x +1)(x^2 +a_4x +1)$   \\
                  {\bf 3} & 4  &    &   & $(x^4 +a_1x^2 +1)(x^4 +a_2x^2 +1)$   \\

  {\bf 4} & 4   &    &   &  $ x^8 +a_1x^4 +1$  \\
  {\bf 5} & 2, 7 &    &   &  $ x^8 -1$  \\
  {\bf 6} &  6  &    &    &  $ x(x^2 +a_1x +1)(x^2 +a_2x +1)(x^2 +a_3x +1)$  \\
  {\bf 7} &  6  &    &  n=2  &  $ x(x^6+a_1x^3+1)$  \\
  {\bf 8} & 7  &    &   &   $(x^2 -1)(x^2 +a_1x +1)(x^2 +a_2x +1)(x^2 +a_3x +1)$ \\
  {\bf 9} & 7  &    &   &  $(x^4 -1)(x^4 +a_1x^2 +1)$  \\
  {\bf 10} &  8,9  &    &   &  $ x(x^6-1)$  \\
  {\bf 11} & 8  &    &   &   $ x(x^2-1)(x^4 +a_1x^2 +1)$ \\
  {\bf 12} & 9  &    &   &  $x(x^2-1)(x^2 +a_1x +1)(x^2 +a_2x +1)$  \\
  {\bf 13} & 12,17  &    &   &  $x^8+14x^4+1$  \\

   \hline
      {\bf 14}  & 1  &   &  & $x^4+ a_1x^3+ a_2x^2+ a_3x +1$  \\
      {\bf 15} & 8  &    &  n=3 & $x(x-1)(x^2 +a_1x +1)$   \\
      {\bf 16} & 8  &    &   & $x(x^3-1)$   \\

     \hline
       {\bf 17} &  4 &   & &  $  (x^2 +a_1x +1)(x^2 +a_2x +1)  $ \\
       {\bf 18} &  4  &    &   &  $x^4 +a_1x^2 +1$  \\
        {\bf 19} & 2,7  &    &   &  $ x^4 -1$ \\
        {\bf 20} &  6 &    &  n=4  &   $x(x^2 +a_1x +1)$ \\   
        {\bf 21} &  7 &    &   &   $(x^2-1)(x^2 +a_1x +1)$ \\ 
        {\bf 22} &  8,9 &    &   &   $x(x^2 -1)$ \\  
        {\bf 23} &  11 &    &   &   $x^4+ 21\sqrt 3 x^2 +1$ \\

       \hline
      {\bf 24}  & 8  &   & n=7 &  $x(x-1)$ \\

 \hline
      \hline
                  
\end{tabular}
      \end{center}
      \end{table}

\clearpage 
\normalsize
 %*********************************************
\section{Concluding remarks}
 
%\subsection{Theta functions of superelliptic curves}   

Finally, we are able to compute for a given genus $g\geq 2$ all full automorphism groups, equations,  of genus $g$ superelliptic curves defined over any algebraically closed field of characteristic different from two. We organize them according to their level $n$.   

These tables are computed for all genus $g \leq 10$ in all characteristics and include for each curve the normal equation of the curve, the automorphism group, the invariants of the corresponding binary forms, $\s$-invariants, the dimension of the corresponding moduli space. Such results will be presented in a continuation of this paper, \cite{super2} where some of the algorithms will be described in more detail.  %Furthermore, all results and some interesting remarks on such computations are planned for \cite{super3}. 

In a current project  we study superelliptic curves defined over $\C$. The half-integer theta characteristics will be computed in each case  and the equation of the corresponding curve in terms of these characteristics, see \cite{beshaj_theta}.

\end{document}

%% file: september.tex
\issueinfo{5}{3}{September}{2011}

\copyrightinfo{2011}{Aulona Press  \emph{(Albanian J. Math.)}}

\def\issn{{\sc \textbf{ISSN: }} 1930-1235; }
\def\issueyear{\textbf{2011}}

\PII{\issn  (\issueyear)}